\documentclass[a4paper,12pt]{article}
\usepackage[T2A]{fontenc}
\usepackage[cp1251]{inputenc}
\usepackage[english,russian]{babel}

\usepackage{graphicx}
\usepackage{amsfonts,amssymb,amscd,amsmath,amsthm}
\usepackage{ccaption}
\captionstyle{\centering}
\captiondelim{. }
\usepackage{latexsym}
\usepackage{lipsum}
\usepackage{indentfirst}
\usepackage{xcolor}
\usepackage[urlbordercolor=white]{hyperref}

\newtheorem{theorem}{Теорема}

\newtheorem{proposition}{Предложение}
\newtheorem{corollary}{Следствие}

\theoremstyle{definition}

\theoremstyle{remark}
\newtheorem{remark}{Замечание}
 
\title{О геометрической медиане треугольной области и других медианоподобных точках}
\author{Петр Панов{$^1$}, Алексей Савватеев$^{2,3}$\vspace{2mm}\\
\small
$^1$Национальный исследовательский университет\vspace{-1.5mm}\\
\small<<Высшая школа экономики>>,  \href{mailto:panovpeter@mail.ru}{panovpeter@mail.ru}\vspace{-1.5mm}\\
\small $^2$Московский физико-технический институт, \href{mailto:hibiny@mail.ru}{hibiny@mail.ru}\vspace{-1.5mm}\\
\small $^3$Центральный экономико-математический институт РАН}

\begin{document}
\maketitle
\begin{abstract}
Геометрическая медиана плоской области --- это точка, минимизирующая среднее расстояние от самой себя до всех точек этой области. Здесь мы выпишем некоторую градиентную систему для вычисления геометрической медианы треугольной области и сформулируем наглядное характеристическое свойство этой медианы. Оно заключается в том, что все три средних расстояния от геометрической медианы до трех сторон границы треугольной области равны между собой. Дальше эти результаты обобщаются на другие виды областей, а также на другие медианоподобные точки.
\end{abstract}
\begin{quote}
\small\textbf{Ключевые слова:} геометрическая медиана, градиентная система, треугольная область.
\end{quote}

Геометрическая медиана является естественным пространственным обобщением статистической медианы одномерной выборки, которая, как известно, минимизирует суммарное расстояние до всех элементов выборки. Именно это минимизирующее свойство положено в основу определения геометрической медианы $m$ конечного набора точек $P_1,\dots P_n$ на плоскости
\begin{equation}\label{eq:nmedian}
	m = \mathop{\mathrm{arg\,min}}\limits_{X\in\mathbb R^2}
	\sum_{i=1}^n{|P_i-X|}.
\end{equation}

С начала прошлого века геометрическая медиана и ее непосредственные обобщения используются в экономической науке в качестве полезного инструмента \cite{Weber1909}. Параллельно продолжают исследоваться математические свойства дискретной медианы и разрабатываются эффективные численные методы для ее нахождения \cite{Wesolowsky1993}. А ближе к концу века интерес смещается в сторону непрерывного случая --- развиваются исследования, связанные с геометрическими медианами кривых и областей \cite{Fekete2005,Zhang2014}. 

В своем изложении мы как раз сосредоточимся на непрерывном случае. Мы начнем с вывода градиентной системы для нахождения геометрической медианы треугольной области (Теорема \ref{TriangleSystem}). Это позволит получить простое и компактное характеристическое свойство геометрической медианы этой области (Предложение \ref{TriangleSystem'}). Дальше эти результаты обобщаются на другие виды областей, а также на другие медианоподобные точки.

Напомним, что по аналогии с дискретным случаем (\ref{eq:nmedian}) геометрическая медиана $m$ области $\Omega\subset \mathbb R^2$ определяется как
\begin{equation}
	m = \mathop{\mathrm{arg\,min}}\limits_{X\in\mathbb R^2}\int_{P\in\Omega}{|P - X|}\,dP,
\end{equation}
где $|P - X|$ --- это обычное евклидово расстояние между точками $P$ и $X$. После введения обозначения
\begin{equation}\label{eq:SigmaOmega}
	\Sigma_\Omega(X) = \int_{P\in\Omega}|P-X|\,dP
\end{equation}
то же самое определение можно записать короче
\begin{equation}
	m = \mathop{\mathrm{arg\,min}}\limits_{X\in\mathbb R^2} \Sigma_\Omega(X).
\end{equation}
Нам понадобится еще одно одно обозначение, пусть $P_1$ и $P_2$ --- некоторые точки на плоскости, тогда
\begin{equation*}
	\Sigma_{P_1P_2}(X) = \int_{P\in P_1P_2}|P-X|\,dP,
\end{equation*}
где интегрирование ведется по отрезку $P_1P_2$. Заметим, что среднее расстояние от точки $X$ до точек отрезка $P_1P_2$ имеет вид $\Sigma_{P_1P_2}(X)/|P_2-P_1|$. Теперь приступим к формулировке первого результата.

\begin{theorem}\label{TriangleSystem}
	Точка $m$ тогда и только будет геометрической медианой треугольной области $\Delta$ с вершинами $P_1,P_2,P_3$, когда выполняется равенство
\begin{equation}\label{eq:TriangleSystem}
	\frac{\Sigma_{P_1P_2}(m)}{|P_2 - P_1|}\,\overrightarrow{P_1P_2} +
	\frac{\Sigma_{P_2P_3}(m)}{|P_3 - P_2|}\,\overrightarrow{P_2P_3} +
	\frac{\Sigma_{P_3P_1}(m)}{|P_1 - P_3|}\,\overrightarrow{P_3P_1} = 0.
\end{equation}
\end{theorem}
{\em Доказательство}. Из определения следует, что медиана $m= m(\Delta)$ --- это критическая точка функции $\Sigma_\Delta$. Возьмем малый вектор $\vec \delta$, $|\vec \delta| = \delta$ и вычислим приращение функции $\Sigma_\Delta$ при смещении аргумента на этот вектор
	\[\delta\Sigma_\Delta(X)=
	\Sigma_\Delta(X + \vec \delta) - \Sigma_\Delta(X).
\]
В критической точке $m$ оно должно иметь порядок $o(\delta)$.

\begin{figure}
\begin{center}
\includegraphics[scale=0.65]{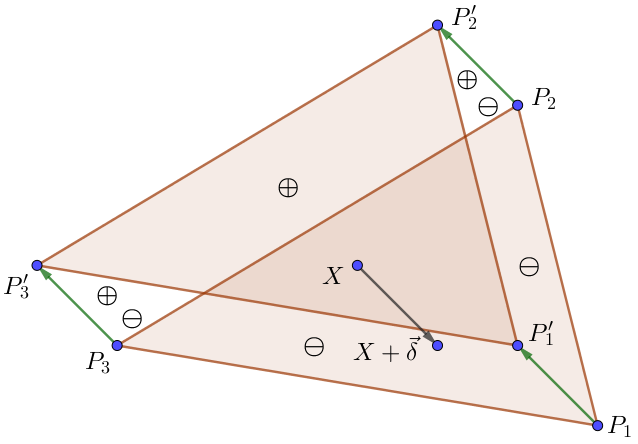}
\vspace{-8pt}
\caption{Треугольник $\Delta'$ сдвинут относительно $\Delta$ на вектор $-\vec \delta$}
\label{fig:BrownTriangle}
\end{center}
\end{figure}

Положим $P'_i = P_i - \vec \delta$ и обозначим через $\Delta'$ сдвинутый треугольник с штрихованными вершинами (рис.~\ref{fig:BrownTriangle}).
Этот рисунок показывает, что приращение 
функции $\Sigma_\Delta$ при смещении аргумента на вектор $\vec \delta$ мы можем записать также и в виде
	\[\delta\Sigma_\Delta(X) =
	\Sigma_{\Delta'}(X) - \Sigma_\Delta(X).
\]

Для параллелограмма с вершинами $P_i,P_{i+1},P'_{i+1},P'_i$ введем обозначение $\pi_i$ и, заодно, введем обозначение~$|\pi_i|$ для его площади. Тогда, как показывает рисунок \ref{fig:BrownTriangle}, это приращение представляется суммой трех слагаемых вида $\Sigma_{\pi_i}(X)$, взятых с подходящими знаками. Причем при $\delta \to 0$ средние значения интегралов $\Sigma_{\pi_i}(X)$ и $\Sigma_{P_iP_{i+1}}(X)$ сближаются друг с другом, а именно,
	\[\frac{\Sigma_{\pi_i}(X)}{|\pi_i|} =
	\frac{\Sigma_{P_iP_{i+1}}(X)}{|P_{i+1}-P_i|} + o(1).
\]
Вопрос о знаке, с которым слагаемое $\Sigma_{\pi_i}(X)$ должно входить в приращение $\delta\Sigma_\Delta(X)$, решается следующим образом. Умножим обе части предыдущего равенства на ориентированную площадь параллелограмма $\pi_i$, а именно, на $-\vec\delta\wedge\overrightarrow{P_iP_i}_{+1}$,
	\[-\frac{\vec\delta\wedge\overrightarrow{P_iP_i}_{+1}}{|\pi_i|}\, \Sigma_{\pi_i}(X) =
	-\vec\delta\wedge\overrightarrow{P_iP_i}_{+1}\,
	\frac{\Sigma_{P_iP_{i+1}}(X)}{|P_{i+1}-P_i|} + o(\delta).
\]
Нетрудно проверить, что дробь, стоящая перед $\Sigma_{\pi_i}(X)$, как раз и является тем самым необходимым знаком. Таким образом,
\begin{multline}\label{eq:deltaSigma}
	\delta\Sigma_\Delta(X) = \\
	= -\vec\delta\wedge\
	\left(
	\frac{\Sigma_{P_1P_2}(X)}{|P_2 - P_1|}\,\overrightarrow{P_1P_2} +
	\frac{\Sigma_{P_2P_3}(X)}{|P_3 - P_2|}\,\overrightarrow{P_2P_3} +
	\frac{\Sigma_{P_3P_1}(X)}{|P_1 - P_3|}\,\overrightarrow{P_3P_1}
	\right) + o(\delta).
\end{multline}
Мы видим, что приращение $\delta\Sigma_\Delta(X)$ в точке $X$  в том и только в том случае будет иметь порядок $o(\delta)$, когда выражение, стоящее в скобках, будет равно нулю. Теорема доказана.
\begin{remark}
На самом деле вектор, расположенный в скобках в выражении (\ref{eq:deltaSigma}) --- это градиент функции $\Sigma_\Delta(X)$, повернутый на $-90^\circ$, так что уравнение (\ref{eq:TriangleSystem}) --- это, действительно, градиентная система.
\end{remark}

Градиентную систему (\ref{eq:TriangleSystem}) можно записать в более компактном и симметричном виде. Из Теоремы \ref{TriangleSystem} вытекает следующее утверждение.

\begin{proposition}[Характеристическое свойство геометрической медианы треугольной области]\label{TriangleSystem'}
Точка $m$ тогда и только тогда будет геометрической медианой треугольной области $\Delta$ с вершинами $P_1,P_2,P_3$, когда
\begin{equation}\label{eq:TriangleSystem'}
	\frac{\Sigma_{P_1P_2}(m)}{|P_2-P_1|} = \frac{\Sigma_{P_2P_3}(m)}{|P_3-P_2|} = \frac{\Sigma_{P_3P_1}(m)}{|P_1-P_3|} .
\end{equation}
Таким образом геометрическая медиана треугольной области --- это точка, для которой три средних расстояния до трех сторон граничного треугольника равны между собой.
\end{proposition}
{\em Доказательство}. В треугольнике одну из сторон выразим через две другие $\overrightarrow{P_3P_1} = -\overrightarrow{P_1P_2} - \overrightarrow{P_2P_3}$, тогда для него равенство (\ref{eq:TriangleSystem}) можно переписать в виде
\begin{equation*}
	\overrightarrow{P_1P_2}\,\left(\frac{\Sigma_{P_1P_2}(m)}{|P_2 - P_1|} - \frac{\Sigma_{P_3P_1}(m)}{|P_1 - P_3|}\right) +
\overrightarrow{P_2P_3}\,\left(\frac{\Sigma_{P_2P_3}(m)}{|P_2 - P_1|} - \frac{\Sigma_{P_3P_1}(m)}{|P_1 - P_3|}\right) = 0.
\end{equation*}
Из-за независимости векторов $P_1P_2$ и $P_2P_3$ следует, что соответствующие им коэффициенты должны быть равны 0, и соотношение (\ref{eq:TriangleSystem'}) в критической точке выполняется. Следствие доказано.
\begin{remark}
Отметим, что все интегралы вида $\Sigma_{PQ}(X)$, присутствующие в градиентных системах (\ref{eq:TriangleSystem}) и (\ref{eq:TriangleSystem'}), вычисляются в конечном виде,  подынтегральная функция тут --- это квадратный корень из квадратного трехчлена. Таким образом, геометрическая медиана это общий ноль двух элементарных функций. Тем не менее,  аналитическая формула конечного вида для нее, по-видимому, неизвестна. Во всяком случае  среди 25\,000 точек, перечисленных в Encyclopedia of Triangle Centers \cite{Kimb} Кларка Кимберлинга, геометрическая медиана треугольной области отсутствует. Заметим, что эта энциклопедия постоянно пополняется и, кроме того, она снабжена проверочными инструментами, позволяющими выяснить, принадлежит ли ей данная точка или нет.
\end{remark}

Точно так же, как и Теорема \ref{TriangleSystem}, доказывается следующий результат.

\begin{figure}[ht]
\centering
\includegraphics[scale=0.65]{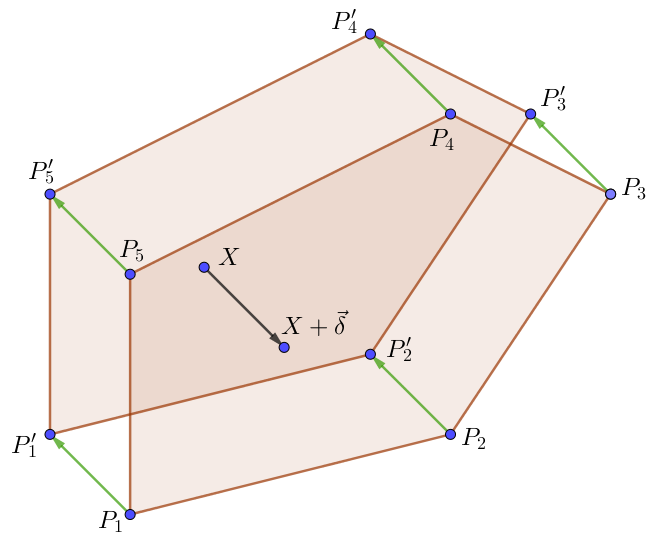}
\caption{Расположение точки $X + \vec\delta$ в многоугольнике $\Pi$ такое же, как у точки $X$ в сдвинутом многоугольнике $\Pi'$}\label{fig:BrownPoly}
\end{figure}

\begin{proposition}\label{PiSystem}
Точка $m$ тогда и только будет геометрической медианой многоугольной области $\Pi$ с вершинами $P_1,\dots,P_n$, когда выполняется равенство
\begin{equation}\label{eq:PiSystem}
	\frac{\Sigma_{P_1P_2}(m)}{|P_2 - P_1|}\, \overrightarrow{P_1P_2}+ \dots +
		\frac{\Sigma_{P_nP_{n+1}}(m)}{|P_{n+1} - P_n|}\,\overrightarrow{P_nP_n}_{+1} = 0.
\end{equation}
\end{proposition}
\noindent
Не приводя доказательства, ограничимся здесь только размещением следующего рисунка~\ref{fig:BrownPoly} --- аналога рисунка \ref{fig:BrownTriangle}.

Для произвольной плоской области $\Omega$, как непосредственное следствие соотношения (\ref{eq:PiSystem}), получаем следующий результат.
\begin{theorem}\label{OmegaSystem}
Пусть задана плоская область $\Omega$ с кусочно гладкой границей. Точка $m$ тогда и только будет ее геометрической медианой, когда выполняется равенство
\begin{equation}\label{eq:OmegaSystem}
	\int_{P\in\partial\Omega}{f(P-m)}\,\overrightarrow{dP} = 0.
\end{equation}
\end{theorem}
\noindent
Чтобы убедиться в правильности этого утверждения, достаточно вписать в кривую $\partial\Omega$ многоугольник с вершинами $P_1,\dots,P_n$, записать для его геометрической медианы равенство (\ref{eq:PiSystem}) и убедиться, что левая часть этого равенства представляет собой интегральную сумму для интеграла~(\ref{eq:OmegaSystem}).

Теперь дополнительно снимем ограничение $n=2$ на размерность и вместо функции $\Sigma_\Omega$, определенной равенством (\ref{eq:SigmaOmega}), будем рассматривать функции более общего вида, определенные равенством
\begin{equation}\label{eq:SimgaOmegaf}
	\Sigma_\Omega^f(X) = \int_{P\in\Omega}{f(P-X)}\,dP,
\end{equation}
где $f$ --- это некоторая функция аргумента $P\in \mathbb R^n$. Здесь имеет место следующий общий результат, который покрывает все предыдущие.
\begin{theorem}\label{th:fOmegaSystem}
Пусть функция $f$ непрерывна и $\Omega$ --- это ограниченная область с кусочно гладкой границей, расположенная в пространстве $\mathbb R^n$. Тогда для функции $\Sigma^f_\Omega$ условие критичности точки $m$ равносильно выполнению следующего равенства
\begin{equation}\label{eq:fromStokes}
	\int_{P\in\partial\Omega}{f(P-m)}\,\overrightarrow{n(P)}\,dP = 0,
\end{equation}
где $\overrightarrow{n(P)}$ --- это единичный вектор внешней нормали к границе области $\Omega$ в точке $P$.
\end{theorem}
\noindent
Доказательство этой теоремы вполне аналогично доказательству Теоремы \ref{OmegaSystem}.

Вернемся к тому, с чего мы начали, а именно, к геометрической медиане треугольной области. Как уже говорилось, в общем случае не удается найти точное решение градиентной системы (\ref{eq:TriangleSystem}) или системы (\ref{eq:TriangleSystem'}). Мы смогли это сделать только в одном тривиальном случае, речь о ``сплюснутом'' треугольнике со сторонами $\alpha,\beta,\gamma$, в котором $\gamma = |\alpha - \beta|$. Имеется ввиду следующее утверждение.
\begin{corollary}
Рассмотрим треугольники со сторонами  $\alpha,\beta,\gamma$, где две большие стороны $\alpha$ и $\beta$ фиксированы, а меньшая $\gamma$ является параметром. Обозначим геометрическую медиану такого треугольника $m(\gamma)$. Пусть теперь $\gamma$ стремится к $|\alpha-\beta|$, тогда $m(\gamma)$ стремится к точке $m$, которая отстоит от общей вершины сторон $\alpha$ и $\beta$  на расстояние~$\sqrt{\alpha\beta/2}$.
\end{corollary}
\noindent
Это утверждение является простым следствием равенства (\ref{eq:TriangleSystem'}), примененного непосредственно к сплюснутому треугольнику со сторонами $\alpha,\beta$ и $|\alpha-\beta|$, для которого оно имеет смысл и для которого в то же самое время само понятие геометрической медианы не вполне определено.

Добавим, что среди сплюснутых треугольников имеются два типа равнобедренных --- это треугольники со сторонами $\alpha,\alpha/2,\alpha/2$ и $\alpha,\alpha,0$. В работе \cite{Panov2018} содержится более точная асимптотическая информация о расположении их геометрических медиан.

\end{document}